\documentclass[conference]{IEEEtran}
\IEEEoverridecommandlockouts
\usepackage{cite}

\usepackage{lipsum}

\ifCLASSINFOpdf
  \usepackage[pdftex]{graphicx}
  \graphicspath{{/}}
  \DeclareGraphicsExtensions{.pdf,.jpeg,.png}
\else
\fi
\usepackage[cmex10]{amsmath}
\usepackage{amssymb}

\usepackage{array}

\usepackage{mdwmath}
\usepackage{mdwtab}
\usepackage[caption=false,font=footnotesize]{subfig}
\usepackage{fixltx2e}

\usepackage{stfloats}

\usepackage{url}

\DeclareMathOperator*{\Ex}{\mbox{$\mathbb{E}$}}
\usepackage{eurosym}
\usepackage{makecell} 

\usepackage{cases}
\renewcommand{\baselinestretch}{0.98}
\begin{document}
%
\title{Medium-term optimization of pumped hydro storage \\with stochastic intrastage subproblems}

\author{\IEEEauthorblockN{Hubert Abgottspon, G\"{o}ran Andersson}
\IEEEauthorblockA{Power Systems Laboratory\\ ETH Zurich, Switzerland\\
Email: \{abgottspon, andersson\}@eeh.ee.ethz.ch}}


%



\maketitle

\begin{abstract}
This paper presents a medium-term self-scheduling optimization of pumped hydro storage power plants with detailed consideration of short-term flexibility. A decomposition of the problem into inter- and intrastage subproblems, where the intrastage problems themselves are formulated as multi-stage stochastic programs, allows the detailed consideration of short-term flexibility. The method is presented together with three alternative approaches, where the short-term flexibility is considered differently: (1) with aggregated peak and off-peak prices, (2) with price duration curves and (3) with deterministic intrastage subproblems. The methods are compared and evaluated in a Monte Carlo operation simulation study. The study is performed on a realistic hydro power plant with consideration of revenue from ancillary services.\end{abstract}



%
\IEEEpeerreviewmaketitle

\section{Introduction}
The goal of a medium-term hydro optimization is to find a seasonal operation strategy. One way to do it is to estimate future revenue, the \emph{profit-to-go}, and to calculate production opportunity costs out of it, the \emph{water values}. 
In Switzerland hydro power plants typically have storage reservoirs which are operated seasonally connected to smaller daily operated reservoirs. 
Therefore future revenues are also influenced by short-term decisions either physically because of e.g.\ empty smaller reservoirs or directly because of the hourly energy market.\\
It is hardly possible to consider a hourly time resolution for a medium-term optimization in a yearly time horizon both computationally as well as because of modeling issues. So aggregations and/or simplifications have to be made where it is difficult not to under- or overvalue short-term flexibility. One important aspect here is how and when information about uncertain variables is disclosed in the model.

\subsection{Proposed model: Stochastic intrastage subproblems: A multihorizon stochastic programming approach}
The proposed modeling approach is based on the obvious observations, that the water management in bigger reservoirs can be considered in a longer time scale than in smaller reservoirs and secondly, that the filling of smaller reservoirs are less important for a revenue estimation. Therefore the idea is that only for seasonal operated reservoirs water values are calculated. The optimization is done for weekly time stages where for each time stage water values are calculated. This multi-stage stochastic program is decomposed into an \emph{interstage} and \emph{intrastage} problem (similar ideas were applied in \cite{Pritchard2005,Abgottspon2012a,Loehndorf2013}): The weekly interstage problem as master problem handles water management in the seasonal reservoirs, whereas the intrastage subproblems ensure hourly water balances in the daily operated reservoirs as well as day-ahead bidding.\\
In contrast to previous works the intrastage problems are not considered deterministically but stochastically. From the modeling point of view this makes sense, since in a weekly perspective hourly water inflows and market prices are not known beforehand. This approach is called, as proposed in another context in \cite{Kaut2013a}, \emph{multihorizon} stochastic programming.\\
Additionally the provision of spinning reserves is modeled as a here-and-now decision within the intrastage problem similar as presented in our previous work \cite{Abgottspon2012a}.\\ 
Summarizing, the proposed model can be described as follows: 
\begin{enumerate}
\item interstage problem (master problem): 
\begin{itemize}
\item weekly water values depending on the filling of the seasonal reservoirs
\item decisions about water release from seasonal reservoirs without information about water inflows and day-ahead prices
\end{itemize}
\item intrastage problem:
\begin{itemize}
\item hourly time steps
\item decisions about provision of spinning reserves, day-ahead bidding and production operation
\item stochastic water inflows are revealed weekly
\item stochastic market prices are revealed daily
\end{itemize}
\end{enumerate}
The consideration of a market for spinning reserves leads to non-concave profit-to-go functions in the master problem (see also Fig. \ref{secondary}). Therefore it is difficult to apply an 
iteration algorithm like stochastic dual dynamic programming. Additionally only a few basins have to be considered in the master problem which eases one of the curses of dimensionality in stochastic dynamic programming \cite{Powell2011}. Further since the goal of this work is to evaluate different modeling methods the solution method is primarily not of interest. Out of these reasons stochastic dynamic programming is applied for the master problem: the problem is decomposed in time and seasonal reservoirs as well as the amount of water discharge from these reservoirs are discretized.\\
The multi-stage stochastic program in the intrastage problem is difficult to decompose since the weekly discharge is given which makes the stages in the intrastage problem depending on each other. The problem is therefore formulated as a deterministic equivalent and a solver based on the simplex method is used to solve it.

\subsection{Evaluation of the proposed method against alternative approaches}
The second goal of this paper is the comparison and evaluation of several methods of aggregations and simplifications for being able to consider short-term flexibility in a medium-term hydro optimization. Apart from the method with stochastic intrastage subproblems the approaches are:
\begin{enumerate}
\item neglecting hourly short-term flexibility with peak and off-peak prices (usual approach)
\item price duration curves (e.g.\ \cite{Pritchard2005,Densing2013})
\item deterministic intrastage subproblems (e.g.\ \cite{Abgottspon2012a,Loehndorf2013})
\end{enumerate}
The four methods are formulated in a way to allow the application of stochastic dynamic programming for the consideration of stochastic inflows and market prices. The methods are evaluated in an operation simulation where the policies suggested by the different optimizations are applied for a number of trial years.\\
The outcome of the evaluation depends on the market structure. Therefore results with- and without the consideration of the revenue out of spinning reserves provision will be shown for a typical hydro power plant with storages.

\subsection{References and Contributions}
Notable references for stochastic programming in the energy sector and stochastic dynamic programming in particular are \cite{Wallace2003a,Masse1946,Bellman1957a,Little1955,Labadie2004a}. The idea of inter- and intrastage problems for hydro power planning were explicitly introduced in \cite{Pritchard2005} and applied in \cite{Loehndorf2013,Abgottspon2012a,Abgottspon2012}, focusing on bidding or operational feasibility respectively.\\
The usual alternative modeling approaches are aggregation and/or bundling of market products (e.g.\ peak and off-peak products) like in \cite{Conejo2008a} and to model different lengths of time stages 
as in \cite{Abgottspon2012,Fleten2011}. Another possibility is shown in \cite{Skjelbred2011} where Lagrangian relaxation was used to incorporate long-term guidelines into short-term and vice versa.\\
Instead of applying approximations to the model another idea is to look for approximate solutions to a detailed model. In \cite{Matos2012} different techniques are reviewed on how to improve the convergence speed of stochastic dual dynamic programming algorithms 
and in \cite{Cerisola2012,Thome2013} water head effects were convexified in order to use such methods. This would also be applicable in this context (see also \cite{Abgottspon2014a}). Depending on the amount of state variables this may be the better approach however modeling issues would still be present.

The contributions of this paper are threefold: First, to the best of our knowledge the application of stochastic intrastage subproblems is presented for the first time. Secondly, several approaches for how to account for short-term flexibility in a medium-term hydro optimization are evaluated and compared. Finally the approaches are extended by considering revenue out of the provision of spinning reserves.

The remainder of the paper is organized as follows: In the next section the model is explained both conceptually and mathematically for all four methods. Afterwards in section III the methods are evaluated and applied in an operation simulation. Some remarks conclude the paper.

\section{Model}
\begin{table} [!t]
\renewcommand{\arraystretch}{1.3}
\caption{Variables}
\label{table}
\centering
\begin{tabular}{ll}
\textbf{Variables} & \textbf{Explanation}\\
\hline
$t \in \mathrm{T} = \{1,...,T\}$ & time stages, $T=52$ [week]\\
$x(t) $ & state and decision variables \\
$\xi_t \in \Xi_t$ & realized data, possibly random\\
$\theta_t(x_{t-1})$ & profit-to-go functions [\euro] \\
$c_t$ & objective function coefficients\\
$\textrm{A}_{t},b_{t}$ & equality constraints (water and financial balances)\\
$\textrm{D}_{t},d_{t}$ & inequality constraints (frequency control reserves)\\
$lb_{t},ub_t$ & lower / upper bounds\\
\hline
$v(t),s(t)$ & filling and spill of reservoirs [m$^3$]\\
$a(t)$ & charges from upstream reservoirs/inflows [m$^3$]\\
$u(t),p(t)$ & generating / pumping [MW] \\
$m(t)$ & position on energy market [MW] \\
$f_u(u),f_p(p)$ & functions of used/produced energy to water flow\\
$q(t)$ & binary: provision of frequency control reserves[0/1]\\
$q^{min}$ & minimum generation [MW]\\
$q^{max}$ & maximum frequency control reserves capacity [MW] \\
$c_t^{q}$ & remuneration for frequency control reserves [\euro/MW] \\
\hline
$c_t^{peak},c_t^{off-peak}$ & aggregated peak and off-peak prices [\euro/MW] \\
$W_t$ & discretized sum of weekly water discharges [m$^3$] \\
\hline
pdc & price duration curve (function) \\
$h^u,h^p$ & hours of generating/pumping per week [h] \\
\hline
$\tau \in \{1,...,168h\}$ & hourly intrastage time steps [h] \\ 
$c_t^{pool}(\tau)$ & hourly day-ahead prices [\euro/MWh] \\ 
$Q_{\xi_t}$ & optimal value of intrastage problem [\euro] \\
$v^{big}(t)$ & weekly filling of seasonal reservoirs [m$^3$]\\
$v^{small}(\tau)$ & hourly filling of daily reservoirs [m$^3$]\\
\hline
$s \in \mathcal{S}$ & scenarios in a scenario tree\\
$\mathcal{A_\tau} \subseteq \mathcal{S}$ & bundle, same past intrastage decisions\\
$\Lambda_\tau$ & set of all bundles $\mathcal{A_\tau}$ in stage $\tau$ \\
$U(\mathcal{A_\tau})$ & set of bundles for ``children'' of $\mathcal{A_\tau}$
\end{tabular}
\end{table}

The overall problem consists of finding a good operating policy for a hydro plant, the \emph{water values} as production opportunity costs. The focus of this paper is the consideration of short-term flexibility where we describe in this section four methods which deal with it differently. Additionally, revenue out of spinning reserves are considered.\\
First model assumptions and limitations are discussed. Then the basic mathematical model is introduced which is then extended for each of the four methods.

\subsection{Model assumptions and limitations}
It is assumed that the power plant is built up out of reservoirs where some are operated in seasonal cycles and others in daily cycles. Further it is assumed, that the sum of the water inflows for one year does not fluctuate much which is valid for power plants in the Alps. Therefore the consideration of a yearly time horizon is sufficient.

As spinning reserves market the one for the provision of secondary frequency control reserves is considered since economically it is the most interesting one in Switzerland. The current market rules require the bidding of symmetric power bands. If the bid is accepted the power band has to be provided for the tendered period of one week. The actual demand is requested automatically and it is assumed, that this request is  symmetric within the tender period. Considered profit out of this market is the remuneration for holding the capacity whereas payments for energy delivery is neglected.\\
When spinning reserves are provided then the turbines have to be continuously running at a certain set point (see Fig. \ref{secondary}). To prevent, that the turbines are operated inefficiently, a minimum generation amount is defined. So the provision of secondary control reserves reduces the production flexibility considerably, which should increase the benefit of a detailed model.

As stochastic variables both water inflows and market prices are considered. The time duration for the main steps is one week, since at the moment the considered spinning reserves are procured weekly in Switzerland. However weekly profit-to-go functions make also sense for a medium-term optimization of the chosen power plant.

\subsection{Mathematical model}
The hydro scheduling problem naturally leads to a multi-stage stochastic program which can be formulated in a dynamic way. Let $\theta_{t}$ be the expected future profit, the \emph{profit-to-go functions}, for the seasonal reservoirs. Then it can be stated:
\begin{align}
\theta_{t}(x_{t-1}) &= \max \ c_t^T x_t + \Ex\limits_{\xi_t \in \Xi_t}\left[\theta_{t+1} \right]   \label{eq::profit-to-go} \\
& \text{subject to:} \nonumber \\
& \quad \textrm{B}_t \cdot x_{t-1} + \textrm{A}_{t} \cdot x_t = b_{t}\nonumber \\
& \quad \textrm{D}_t \cdot x_t \le d_{t}\nonumber \\
& \quad lb_t \le x_t \le ub_t \nonumber \\
& \quad x_t \in \mathbb{R}^n,\{0,1\} \nonumber 
\end{align}
$x_{t}$ specifies the state and decision variables at time stage $t$, i.e. the fillings of the reservoirs as well as production and bidding decisions. $\xi_t := (c_t,\textrm{A}_t,\textrm{B}_t,b_t,\textrm{D}_t,d_t)$ defines the data vector where $c_t,b_t$ are random and not known in advanced. These are market prices as well as water inflows.\\
$\Ex[..]$ denotes the expected value over sampled random data $\xi_{t} \in \Xi_t$ and is maximized in order to find the profit-to-go function. Note, that the stochastic data process $\xi_1,...,\xi_T$ is Markovian, so the profit-to-go function $\theta_t$ depends only on $\xi_t$ and not on the whole past process $\xi_1,...,\xi_t$.\\ 
The stochastic program is subject to equality constraints, defined by $\textrm{B}_t,\textrm{A}_t,b_t$. In more detail, these constraints ensure correct water and financial balance. The water balance is modeled as follows:
\begin{equation}
v_t = v_{t-1} -s_t - f_u(u_t) + f_p(p_t) + a_t -q_t \cdot f_u(q^{min} + q^{max})
\label{eq::water_balance}
\end{equation}
To keep notations simple $a(t)$ denotes both water inflows as well as charge from upstream reservoirs.\\
The market position $m(t)$ is maximized in the objective function multiplied with the market prices. It is defined as:
\begin{equation}
m_t = u_t - p_t + q_t \cdot(q^{min} + q^{max})
\label{eq::financial_balance}
\end{equation}
The provision of secondary control reserves influences the operation of the turbines. This is can be modeled as inequality constraints ($\textrm{D}_t,d_t$):
\begin{equation}
 q_t \cdot (q^{min} + q^{max}) \le u_t \le \overline{u} - q_t \cdot q^{max}
\label{eq::secondary}
\end{equation} 
Note that the provision of spinning reserves is approximated by taking into account either no provision or the maximum quantity for each turbine of the power plant. The remuneration for holding a capacity $c_t^q$ is assumed to be known beforehand, it is estimated as the minimum amount one can expect to get accepted. This remuneration is part of the objective function.\\
Finally the lower and upper bounds are:
\begin{align}
&0             \le  \ v_t \le \overline{v} , v_t \in \mathbb{R}^v \quad , \quad 
0			  \le  \ s_t \le \overline{s} , s_t \in \mathbb{R}^v \label{eq::bounds} \\
&0			  \le  \ u_t \le \overline{u} , u_t \in \mathbb{R}^u \quad , \quad
0             \le  \ p_t \le \overline{p} , p_t \in \mathbb{R}^p \nonumber \\
&0			  \le  \ q_t \le 1			 , q_t \in \mathbb{Z}^u \quad , \quad 
\underline{m} \le  \ m_t \le \overline{m} , m_t \in \mathbb{R} \nonumber 
\end{align}
The position $m_t$ is bounded to some upper and lower values to prevent extreme positions and therefore market speculations are prevented.\\
In the following paragraphs this basic model \eqref{eq::profit-to-go} is adapted for the four mentioned methods.

\subsection{Method 1: Neglecting hourly flexibility with weekly peak and off-peak prices}
\begin{figure} [!t]
\centering
\includegraphics[width=1\columnwidth]{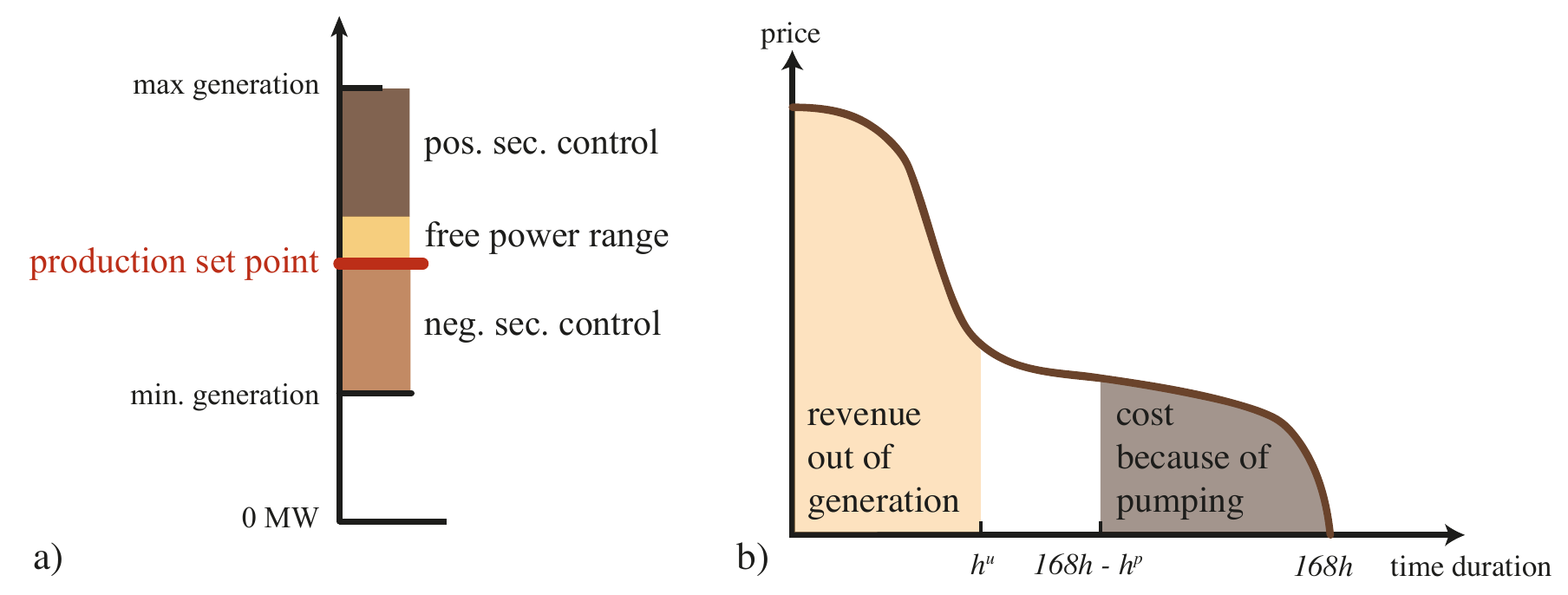}
\caption{a) Generation of a turbine with provision of secondary frequency control reserves. Production flexibility is reduced considerably and the substantial minimum of generation results to non-concave profit-to-go functions. b) Schematic example of a price duration curve. 
Revenue out of generation as well as costs because of pumping are shown. Note that since the overall water discharge is fixed, with more pumping more generation is possible.}
\label{secondary}
\end{figure}
The first method neglects hourly flexibility. Water inflows and market prices are estimated as expected values over the whole next period. Two weekly prices are assumed, peak prices $c_t^{peak}$ when energy generation is done and off-peak prices $c_t^{off-peak}$ for pumping. Smaller basins are disregarded and turbines and pumps are aggregated accordingly (see also Fig. \ref{hydro_plant} b)). This result to problem \eqref{eq::profit-to-go}, where $x_t$ consists of one value for each of the variables $v_t,s_t,a_t$ per aggregated reservoir and one value for $u_t,p_t,q_t$ for each turbine and pump respectively:
\begin{align}
\theta_{t}(x_{t-1}) &= \max_{W_t} \ \Ex\limits_{\xi_t \in \Xi_t} [ c_t^{peak}u_t - c_t^{off-peak}p_t \ ... \\
& \quad  + q_t \cdot q^{max} \cdot c_t^q \cdot 168h + \theta_{t+1}(x_t) ]   \nonumber  
\end{align}
\quad subject to:

\quad \quad contraint \eqref{eq::secondary} and bounds \eqref{eq::bounds} as well as:
\begin{align}
&v_t = v_{t-1} -W_t \label{eq::simple_2} \\
&f_u(u_t) - f_p(p_t) + s_t - a_t= W_t \label{eq::simple_3}
\end{align}
Note that by discretizing the weekly water discharge from the reservoirs $W_t$ it is possible to apply the stochastic dynamic programming scheme. Note also, that random data consists of peak and off-peak prices as well as water inflows $a_t$. \\
The advantage with this formulation is the small computational burden although stochasticity is taken into account. For every time stage $t$ and $W_t$ there is only one single constraint of \eqref{eq::secondary}, \eqref{eq::simple_2} and \eqref{eq::simple_3} (for each scenario).

\subsection{Method 2: Price duration curves}
The price duration curve (example in Fig. \ref{secondary} b)) is constructed out of the proportion of hourly prices below a certain price for some time duration. Since the revenue is something like \emph{price x quantity} of sold energy it can be estimated by integration of the price duration curve.\\
In \cite{Pritchard2005} such curves are multiplied by quantity-price offers and integrated in respect to prices. Here another approach is followed, where the sum of the water discharge for the next week is discretized. Then for a given water discharge an optimization problem is formulated in order to find the time duration of pumping $h_t^p$ and generating $h_t^u$. For given price duration curves the expected short-term profit can be derived.

It is assumed, that the power plant either produces or pumps fully or not for each hour which is a valid approximation in this context. Random data involve again prices and water inflows. To estimate random price duration curves the hourly price process itself is sampled and the price duration curve is constructed out of it. The problem can now be formulated as follows:
\begin{align}
\theta_{t}(x_{t-1}) &= \max_{W_t} \Ex\limits_{\xi_t \in \Xi_t} \Bigg[ \ u_t \cdot \int_{0}^{h_t^u} \text{pdc}_t(\tau) d\tau \ ... \label{eq::pdc} \\
& \quad - \overline{p} \cdot \int_{168h - h_t^p}^{168h} \text{pdc}_t(\tau) d\tau \ ... \nonumber \\
& \quad +  q_t \cdot \Big( (q^{max}+q^{min})\int_{0}^{168h} \text{pdc}_t(\tau) d\tau \ ... \nonumber \\
& \quad + q^{max}c_t^q \cdot 168h \Big) +  \theta_{t+1}(x_t) \Bigg]   \nonumber  
\end{align}
subject to:

\quad \quad bounds \eqref{eq::bounds} as well as:
\begin{align}
&v_t = v_{t-1} - W_t \nonumber \\
&u_t = \bar{u} - q_t\cdot (q^{max} + q^{min}) \nonumber \\
 &h_t^u \sum f_u(\max(u_t)) - h_t^p \sum f_p(\overline{p}) \ ... \nonumber \\
 & \quad + q_t f_u(q^{max}+q^{min}) \cdot 168h + s_t - a_t= W_t \nonumber
\end{align}
The problem turns out to be challenging to solve. Therefore the price duration curves are assumed to be pice-wise linear which approximates the problem to a quadratic mixed-integer problem.

From the modeling point of view there are several approximations with this formulation. The most severe is that similar to the first method timing is not respected at all, which means it is not considered when and in which order the decisions are taken within a week.\footnote{For instances it may happen, that high market prices occur all at the beginning of a week where the reservoirs may be empty and generating not possible which is not taken care of with the consideration of price duration curves.}\\
The advantage with this formulation is the consideration of a reasonable representation of the opportunities in the hourly day-ahead market.

\subsection{Method 3: Deterministic intrastage subproblems}
The idea of intrastage subproblems is already explained in the introduction. For the third method these subproblems are modeled deterministically.\\
Mathematically the multi-stage stochastic program with intrastage subproblems can be formulated in a similar way to \eqref{eq::profit-to-go}:\begin{equation}
\theta_{t}(v_{t-1}^{big}) = \max_{W_t} \ \Ex\limits_{\xi_t \in \Xi_t}\left[ Q_{\xi_t,W_t} + \theta_{t+1}(v_{t}^{big}) \right]
\label{eq::det-intra}
\end{equation}
where:
\begin{align}
&Q_{\xi_t,W_t}(v_{t-1}^{big}) = \max_{u_\tau,p_\tau,s_\tau,q_t} \ (c_t^{pool})^T m_t(\tau) \ ...  \label{eq::det_2nd}\\
& \quad  \quad  + c_t^q \cdot q_t^T q^{max} \cdot 168h \nonumber
\end{align}
\quad subject to:
\begin{align}
&\text{1. water balances:} \nonumber \\
&\quad v_t^{big} = v_{t-1}^{big} -W_t \nonumber \\ 
&\quad v_\tau^{small} = v_{\tau-1}^{small} -s_\tau - f_u(u_\tau) + f_p(p_\tau) + a_\tau  \ , \ \forall \tau \nonumber \\ 
&\quad \sum_\tau \left[f_u(u_\tau) - f_p(p_\tau) - a_\tau \right] + s_t = W_t \nonumber \\
&\text{2. financial balance:} \nonumber \\
&\quad m_\tau = u_\tau - p_\tau \ , \ \forall \tau \nonumber \\
&\text{3. secondary control provision:} \nonumber \\
&\quad q_t \cdot (q^{min} + q^{max}) \le u_\tau \le \overline{u} - q_t \cdot q^{max} \nonumber \\
&\text{4. bounds similar to \eqref{eq::bounds}:} \nonumber \\
&\quad  lb_t \le u_\tau,p_\tau,s_\tau,q_t,v_t^{big},v_\tau^{small},m_\tau \le ub_t \nonumber \\
&\quad  u(\tau),p(\tau),s(\tau),v^{small}(\tau) \in \mathbb{R}^{n\text{x}\tau} ,\ m(\tau) \in \mathbb{R}^{\tau} \nonumber \\
&\quad  v_t^{big} \in \mathbb{R}^n , \ q_t \in \{ 0,1 \}^n \nonumber
\end{align}
$Q_{\xi_t,W_t}$ is the optimal value of the deterministic intrastage subproblem. It is a function of the former state, realized random data $\xi_t$ as well as the discretized water release $W_t$. \\
The purpose of the problem $Q$ is to estimate the intrastage profit in a realistic way, by hourly deploying $W_t$ most optimally within the week. It is formulated as a two-stage stochastic program. In the first stage the amount of secondary control reserves to bid is decided. This is done for each turbine, which is qualified for this provision. Afterwards $\xi_t$ is disclosed, so the water inflows and prices for the whole week become known. Then actual hourly production decisions take place. As a consequence $Q_{\xi_t,W_t}$ is a deterministic, linear maximization problem, with binary variables.

Note, that the operation of the power plant is considered in hourly resolution as opposed to the first and second method. Approximations made are first that the random data are assumed to be known one week in advanced, which is an optimistic view. Further the fillings of the daily reservoirs $v^{small}$ are neglected in calculation of the profit-to-go functions as well as the water balance is not respected between consecutive weeks. So this means that the fillings of the daily reservoirs are zero at the beginning and end of each week.

\subsection{Method 4: Stochastic intrastage problems}
\begin{figure} [!t]
\centering
\includegraphics[width=0.75\columnwidth]{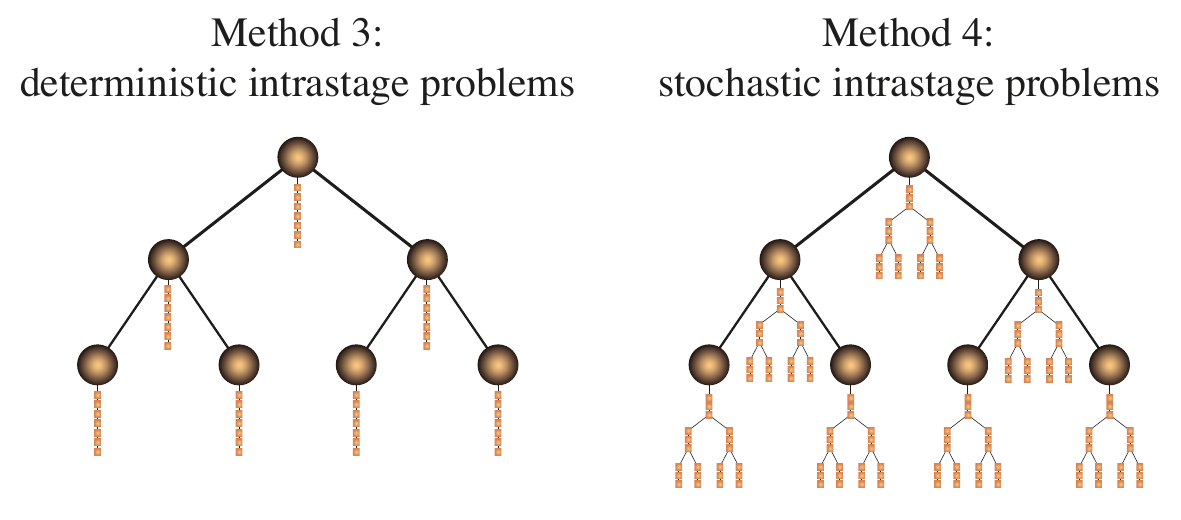}
\caption{Decision trees with deterministic and stochastic intrastage subproblems. Whereas random data in the third method is revealed once for each intrastage problem in the forth method it is revealed daily. Note that interstage decisions do not depend directly on intrastage decisions of previous stages.}
\label{multi-horizon}
\end{figure}
Now the model from the previous method 3 is extended by considering stochastic instead of deterministic intrastage subproblems which is one of the novel contributions of this paper. The idea is depicted in Fig. \ref{multi-horizon}. Whereas in the third method the random data is disclosed at the beginning of a week the day ahead prices are now revealed daily for one day. This seems more realistic, that only market prices for the day ahead are known in advance. The water inflows however are still revealed weekly out of two reasons: to give a hint about the modeling flexibility and to keep computational burden low.\\
Note that the stochastic intrastage problems cannot be formulated in a dynamic way since the sum of the weekly discharge $W_t$ has to be fulfilled for each scenario. So it is formulated and solved as a deterministic equivalent linear program.

The mathematical formulation for stochastic intrastage subproblems changes in respect to how random data is disclosed if compared to \eqref{eq::det-intra} and \eqref{eq::det_2nd}. To keep notation simple the subproblems are described with the help of scenario trees. A scenario is one possible realization path of the random data. Let the set of all scenarios $s$ be $\mathcal{S}$ and consider a \emph{bundle} $\mathcal{A}_\tau\subseteq \mathcal{S}$ a subset of $\mathcal{S}$ with the same intrastage decisions up to some stage $\tau$. Finally let $\Lambda_\tau$ be the set of all bundles in a stage $\tau$ and therefore $\mathcal{A}_\tau \in \Lambda_\tau$. Further let the set of bundles $U(\mathcal{A}_\tau)$ be:\footnote{As an example consider Fig. \ref{multi-horizon}: For each stochastic intrastage subproblem the cardinality of $\mathcal{S}$ is 4, which means there are 4 different scenarios. $\Lambda_4$ is consisted of two bundles: $\Lambda_4 = \Lambda_5 = \Lambda_6 = \{\{1,2\},\{3,4\}\}$. Consider now the bundle $\mathcal{A}_6 = \{1,2\}$. $U(\mathcal{A}_6)$ then specifies the children of $\mathcal{A}_6$, the set of bundles $U(\mathcal{A}_6) = \{\{1\},\{2\}\}$.}
\begin{equation*}
U(\mathcal{A}_\tau) = \left\{\mathcal{B} \in \Lambda_{\tau+1}\, | \, \mathcal{B} \subseteq \mathcal{A}_\tau \right\}
\label{eq::children}
\end{equation*}
Now the formulation of the stochastic intrastage subproblems looks similar to the deterministic one, however there are now not only on variable per stage $\tau$ but for each bundle $\mathcal{A}_\tau$. The problem can be written as follows:
Aendere eventuell zwei Sachen hier: Erstenes $\theta_{t+1}(v_{t-1}^{big}-W_t)$ und zweitens $\Ex\limits_{\mathcal{A}_\tau \in \Lambda_\tau}\left[c_t^{pool})^T \cdot m_{\mathcal{A}_\tau} \right]$. Und dann schau noch, ob pumping / generating vs. used/produced energy indexes korrekt sind.

\begin{equation}
\theta_{t}(v_{t-1}^{big}) = \max_{W_t} \ \Ex\limits_{\xi_\tau \in \Xi_\tau}\left[ Q_{W_t} + \theta_{t+1}(v_{t}^{big}) \right]
\label{eq::stoch-intra}
\end{equation}
where:
\begin{align}
&Q_{W_t}(v_{t-1}^{big}) = \max_{u_\tau,p_\tau,s_\tau,q_t} \ (c_t^{pool})^T \Ex\limits_{\mathcal{A}_\tau \in \Lambda_\tau}\left[m_{\mathcal{A}_\tau} \right] \ ...  \label{eq::stoch_2nd}\\
& \quad  \quad  + c_t^q \cdot q_t^T q^{max} \cdot 168h \nonumber
\end{align}
\quad subject to:
\begin{align}
&\text{1. water balances:} \nonumber \\
&\quad v_t^{big} = v_{t-1}^{big} -W_t \nonumber \\ 
&\quad v_{\mathcal{B_\tau}}^{small} = v_{\mathcal{A}_{\tau-1}}^{small} -s_{\mathcal{B}_\tau} - f_u(u_{\mathcal{B}_\tau}) + f_p(p_{\mathcal{B}_\tau}) + a_{\mathcal{B}_\tau} \ , \ ... \nonumber \\
&\quad \quad \quad \forall \mathcal{A}_{\tau-1} \in \Lambda_{\tau-1}, \ \forall \mathcal{B}_\tau \in U(\mathcal{A}_{\tau-1}) , \ \forall \tau \nonumber \\ 
&\quad \sum_\tau \left[f_u(u_\tau) - f_p(p_\tau) - a_\tau \right] + s_t = W_t \nonumber \\
&\text{2. financial balance:} \nonumber \\
&\quad m_{\mathcal{A}_\tau} = u_{\mathcal{A}_\tau} - p_{\mathcal{A}_\tau} \ , \   \forall \mathcal{A}_\tau \in \Lambda_{\tau} , \ \forall \tau \nonumber \\
&\text{3. secondary control provision:} \nonumber \\
&\quad q_t \cdot (q^{min} + q^{max}) \le u_{\tau,\mathcal{A}_\tau} \le \overline{u} - q_t \cdot q^{max} \ , \ ...\nonumber \\
&\quad \quad \quad  \forall \mathcal{A}_\tau \in \Lambda_{\tau} , \ \forall \tau  \nonumber \\
&\text{4. bounds similar to \eqref{eq::bounds}:} \nonumber \\
&\quad  lb_t \le u_{\mathcal{A}_\tau},p_{\mathcal{A}_\tau},s_{\mathcal{A}_\tau},q_t,v_t^{big},v_{\mathcal{A}_\tau}^{small},m_{\mathcal{A}_\tau} \le ub_t \nonumber \\
&\quad  u(\mathcal{A}_\tau),p(\mathcal{A}_\tau),s(\mathcal{A}_\tau),v^{small}(\mathcal{A}_\tau) \in \mathbb{R}^{n\text{x}\sum_\tau |\Lambda_\tau |} \nonumber \\
&\quad  m(\mathcal{A}_\tau) \in \mathbb{R}^{\sum_\tau |\Lambda_\tau |} , \ v_t^{big} \in \mathbb{R}^n , \ q_t \in \{ 0,1 \}^n \nonumber
\end{align}
\begin{figure} [!t]
\centering
\includegraphics[width=0.89\columnwidth]{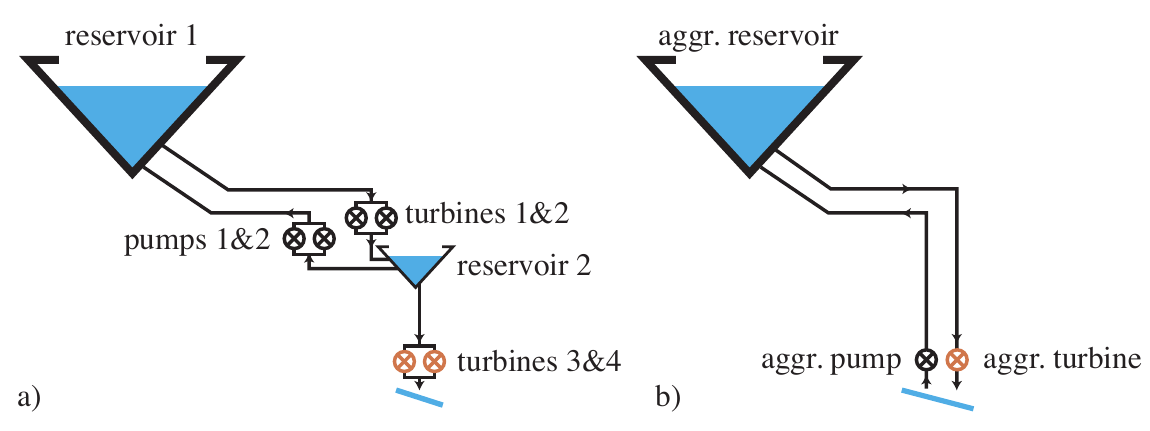}
\caption{a) Schematic overview of the considered hydro power plant. Reservoir 1 is the seasonal storage whereas reservoir 2 is the daily one. Turbines 3 and 4 are qualified to deliver secondary frequency control reserves. b) Aggregated hydro power plant for the first and second method. Note, that also water inflows are aggregated and that an infinite lower basin for pumping is assumed. The aggregated turbine is able to provide frequency control reserves.}
\label{hydro_plant}
\end{figure}
Note, that the size of the variable vectors in the intrastage problem depends on the sum of the number of bundles for each time step $\sum_\tau |\Lambda_\tau |$. For instances consider two day-ahead price scenarios per day. Each intrastage vector (for each reservoir etc.) then have $24 \cdot (2^1+2^2+...+2^7) = 6096$ entries (for comparison the deterministic formulation would require only $24 \cdot 7 = 168$ entries).\\
Compared with the method 3 with deterministic intrastage subproblems the method 4 is much more realistic from the modeling point of view. Computationally the same amount of subproblems have to be solved, however the size of these subproblems differ.\footnote{For example for two daily market price scenarios (which results to $2^7$=128 weekly scenarios) one stochastically formulated intrastage subproblem was constructed and solved in 0.41\,s whereas the deterministic variant required only 0.12\,s.}

\section{Evaluation}
There are no standardized hydro power plant models available for optimization studies as it is the case e.g.\ for electricity grid analysis. The outcome of the evaluation is however depending on the considered power plant. It seems obvious, that the more complicated structure and the more hourly dynamics are present in the model, the better the proposed method with stochastic intrastage subproblems should perform. Chosen was therefore a typical Swiss hydro power plant (Fig. \ref{hydro_plant}) which is not overly complicated however still consists of two different kinds of reservoirs, pumps and turbines as well as is qualified to provide secondary frequency control reserves. Another option would have been to consider several different kind of power plants, but this was beyond the scope of this paper.

The advantages and disadvantages from the modeling point of view are already explained in the previous section. Compared is now the computational burden of the methods, their proposed water values as well as a simulation study where the water values are applied for several samples of water inflows and market prices.

\subsubsection{Computational burden}
\begin{figure} [!t]
\centering
\includegraphics[width=1\columnwidth]{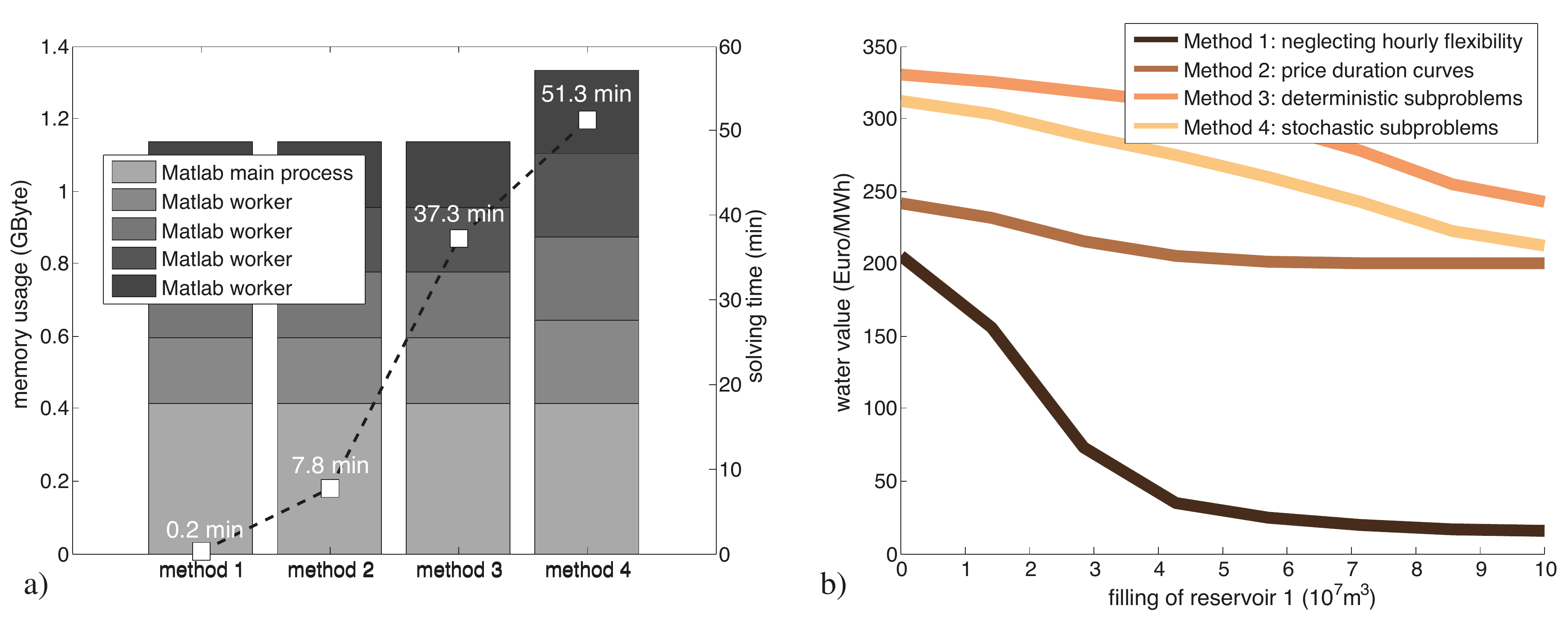}
\caption{a) Time duration needed for solving the different optimization methods and memory requirements of the optimizations. b) Water values (gradients of the profit-to-go functions) for the different methods for the first time stage. They depend on reservoir filling and weekly time stages. Note that the first two methods relatively undervalue short-term revenue whereas the third one overvalues it.}
\label{computational_burden}
\end{figure}
The optimizations were done on a standard computer with 4 physical processor cores. CPLEX dual simplex solvers was used for the linear and quadratic programs whereas for the mixed-integer problem a branch-and-cut algorithm was used.\\
In Fig. \ref{computational_burden} the time duration as well as the memory requirements needed for the different methods are depicted. The first method is solved in 15\,s which is 30 to 200 times faster than the other methods, which would be clearly a big advantage in daily use of such optimizations. It should be noted that it is only method four that has higher memory requirements. Memory usage of this method (as well as solving time) will further increase exponentially by increasing the amount of intra time stages, stochasticity or power plant complexity (state variables).

\subsubsection{Water values comparison}
The results of the optimization methods, the water values, are shown and compared in Fig. \ref{computational_burden} b) for the first time stage. Notable is, that the more short-term dynamics are considered the higher the water value is. However because the third method assumes perfect weekly knowledge it overvalues the water value. Therefore method 4 should give a more realistic estimation. The water values for methods 3 and 4 are similar not only for the first time stage but also for the others with roughly $80\, \%$ of them having a difference of less than $10\,\%$.

\subsubsection{Simulation study}
\begin{figure} 
\centering
\includegraphics[width=1\columnwidth]{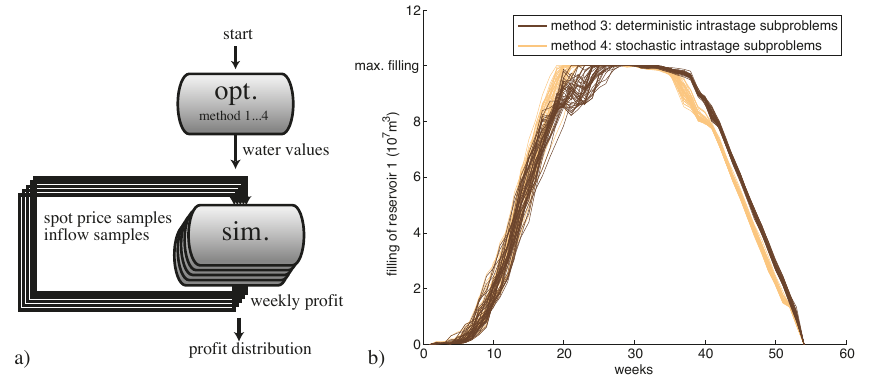}
\caption{a) Flowchart of the Monte Carlo operation simulation. 
b) Filling of seasonal reservoir 1 for each of the samples for water values from method 3 and 4. Note that although the water values are similar the actual operation varies more.}
\label{flow}
\end{figure}
\begin{table*} [!t]
\renewcommand{\arraystretch}{1.1}
\caption{Comparison of optimization methods without / with provision of secondary frequency control reserves}
\label{comparison}
\centering
\begin{tabular}{r|cccc}
&Method 1 & Method 2 & Method 3 & Method 4 \\
&weekly peak \& off-peak prices & price duration curves & deterministic intrastage & stochastic intrastage \\
\hline
 expected profit &34.24 / 34.99 & 29.13 / 30.44 & 34.80 / 35.85 &33.40 / 39.14 \\
 rel. standard deviation  &2.19\% / 2.53\% & 1.95\% / 2.47\% & 2.36\% / 3.99\% &2.57\% / 4.10\% \\
 CV$@$R$_{10\%}$  &32.66 / 33.27 & 27.90 / 29.03 & 33.27 / 33.47 &31.82 / 36.28 \\
\end{tabular}
\end{table*}
A comparison of the water values does not answer the question how well an application of it would perform. A Monte Carlo simulation study shall estimate this performance of the different methods (Fig. \ref{flow} a)). In the simulation part, the hydro power plant operation is mimicked over one year. So for each week realistic hydro power plant operation is simulated, based on the estimated weekly water values. The procedure is as follows:
\begin{enumerate}
\item sampling of correlated water inflows and market prices
\item offering of secondary frequency control reserves
\item hourly production decision heuristic
\end{enumerate}
Because of lack of sufficient amount of historical data distributions are estimated and water inflows and market prices are sampled out of it. Then it is decided, if secondary frequency control reserves are offered for the next week or not which is modeled as a mixed-integer problem. After that a heuristic performs hourly production decisions, which simulates what an operator would do in practice: First frequency control reserves obligations are fulfilled and then energy is generated or used for pumping depending on a comparison of the filling depended water values and market prices. This procedure is repeated for every week and for 100 samples in a receding horizon. Outcome of the simulation is a profit distribution.\\
Fig. \ref{flow} b) shows for method 3 and 4 the resulting seasonal reservoir filling for  all samples. The reservoirs maximum filling is exploited with both methods (without spilling). One could argue, that the application of water values from method 4 leads to a more conservative strategy i.e. releasing water earlier.

Table \ref{comparison} shows the expected profit, the relative standard deviation as well as the mean profit for the 10\% worst scenarios (CV$@$R$_{10\%}$). The values are shown with and without consideration of provision of control reserves.\\
Method 1 leads to astonishingly good results. However the performance evaluation was based on market data, where peak and off-peak price periods were clearly present which will or already is not anymore the case.\\
Method 2 performs worse than expected. An explanation for it could be, that although the price process is considered in a detailed way the power plant itself is simplified considerably, which leads to using non-existing resources more efficiently. This may result to policies which are less effective.\\
Method 3 outperforms method 4 for optimizations without the consideration of secondary frequency control provision. Interesting is also the increased robustness if compared with method 1: the CV$@$R$_{10\%}$ is considerably higher whereas the relative standard deviation as another risk measure would indicate slightly more risk.\\
Finally the proposed method 4 outperforms the other methods only if secondary control reserves provision is considered. But in this case the increase of both expected profit and CV$@$R$_{10\%}$ is substantially by around 10\,\%.

\section{Conclusions}
This paper presented four aggregation methods for a medium-term self-scheduling of hydro power plants. The methods were: (1) aggregated peak and off-peak prices, (2) price duration curves, (3) deterministic intrastage subproblems and (4) stochastic intrastage subproblems. Contributions were first the application of stochastic intrastage subproblems to the hydro power planning, second the comparison and evaluation of the different methods and finally the extension of the methods by considering revenue out of provision of spinning reserves.\\
The evaluation presented the different computational burden as well as the quality of the proposed optimization outcomes, where a Monte Carlo operation simulation study was performed. The results suggest, that the consideration of stochastic intrastage subproblems makes only sense if the market structure is sufficiently complex as it is the case for the market for secondary frequency control provision. The results also indicate that another reason could be more complex hydro plant structures, whereas the evaluation on such a plant is left for future work.

\bibliographystyle{IEEEtran}
\bibliography{literature.bib}

\end{document}